\documentclass[MSNbibl,number,citesort,secthm,seceqn,dvips]{arxbj}
\usepackage{mathbh}

\aid{0}
\volume{18}
\issue{3}
\pubyear{2012}
\firstpage{764}
\lastpage{782}
\doi{10.3150/11-BEJ387}

\makeatletter
\def\bsuffix #1{#1}

\newcommand{\eqref}[1]{(\ref{#1})}
\newcommand{\dto}{\rightsquigarrow}
\newcommand{\pto}{\mathop{\rightarrow}^{p}}
\newcommand{\CC}{\mathbb{C}}
\newcommand{\reals}{\mathbb{R}}
\newcommand{\eps}{\varepsilon}
\newcommand{\norm}[1]{\|{#1}\|}
\newcommand{\ind}{\mathbh{1}}
\newcommand{\cov}{\operatorname{cov}}
\newcommand{\corr}{\operatorname{corr}}
\newcommand{\diff}{\mathrm{d}}
\newcommand{\I}{\mathrm{I}}
\newcommand{\II}{\mathrm{II}}
\newcommand{\III}{\mathrm{III}}

\renewcommand{\ge}{\geq}
\renewcommand{\le}{\leq}
\newremark{example}[proposition]{Example}
\newproclaim{condition}[proposition]{Condition}
\newtheorem{proposition}{Proposition}[section]
\newtheorem{lemma}[proposition]{Lemma}

\makeatother

\begin{document}
\begin{frontmatter}

\title{Asymptotics of empirical copula processes under non-restrictive smoothness assumptions}
\runtitle{Empirical copula processes}

\begin{aug}
\author{\fnms{Johan} \snm{Segers}\corref{}\ead[label=e1]{Johan.Segers@uclouvain.be}}

\address{Universit\'e catholique de Louvain, Institut de statistique, biostatistique et
sciences actuarielles, Voie du Roman Pays 20, B-1348 Louvain-la-Neuve, Belgium. \printead{e1}}
\end{aug}

\received{\smonth{12} \syear{2010}}
\revised{\smonth{4} \syear{2011}}

%
\begin{abstract}
Weak convergence of the empirical copula process is shown to hold under
the assumption that the first-order partial derivatives of the copula
exist and are continuous on certain subsets of the unit hypercube. The
assumption is non-restrictive in the sense that it is needed anyway to
ensure that the candidate limiting process exists and has continuous
trajectories. In addition, resampling methods based on the multiplier
central limit theorem, which require consistent estimation of the
first-order derivatives, continue to be valid. Under certain growth
conditions on the second-order partial derivatives that allow for
explosive behavior near the boundaries, the almost sure rate in Stute's
representation of the empirical copula process can be recovered. The
conditions are verified, for instance, in the case of the Gaussian
copula with full-rank correlation matrix, many Archimedean copulas, and
many extreme-value copulas.
\end{abstract}

%

\end{frontmatter}

\section{Introduction}

A flexible and versatile way to model dependence is via copulas. A
fundamental tool for inference is the empirical copula, which basically
is equal to the empirical distribution function of the sample of
multivariate ranks, rescaled to the unit interval. The asymptotic
behavior of the empirical copula process was studied in, amongst
others, Stute \cite{stute1984}, G{\"a}nssler and Stute
\cite{gaensslerstute1987},
Chapter~5, van~der Vaart and Wellner \cite{vandervaartwellner1996}, page~389,
Tsukahara \cite{tsukahara2000,tsukahara2005},
Fermanian \textit{et al.} \cite{fermanianradulovicwegkamp2004},
Ghoudi and R{\'e}millard \cite{ghoudiremillard2004}, and
van~der Vaart and Wellner~\cite{vandervaartwellner2007}. Weak
convergence is shown typically for copulas that are continuously
differentiable on the closed hypercube, and rates of convergence of
certain remainder terms have been established for copulas that are
twice continuously differentiable on the closed hypercube.
Unfortunately, for many (even most) popular copula families, even the
first-order partial derivatives of the copula fail to be continuous at
some boundary points of the hypercube.

\begin{example}[(Tail dependence)]\label{extaildep}
Let $C$ be a bivariate copula with first-order partial derivatives $\dot
{C}_1$ and $\dot{C}_2$ and positive lower tail dependence coefficient
$\lambda= \lim_{u \downarrow0} C(u, u)/u > 0$. On the one hand, $\dot
{C}_1(u, 0) = 0$ for all $u \in[0, 1]$ by the fact that $C(u, 0) = 0$
for all $u \in[0, 1]$. On the other hand, $\dot{C}_1(0, v) = \lim_{u
\downarrow0} C(u, v)/u \ge\lambda> 0$ for all $v \in(0,1]$. It
follows that $\dot{C}_1$ cannot be continuous at the point $(0, 0)$;
similarly for~$\dot{C}_2$. For copulas with a positive upper tail
dependence coefficient, the first-order partial derivatives cannot be
continuous at the point $(1, 1)$.
\end{example}

Likewise, for the Gaussian copula with non-zero correlation parameter
$\rho$, the first-order partial derivatives fail to be continuous at
the points $(0, 0)$ and $(1, 1)$ if $\rho> 0$ and at the points $(0,
1)$ and $(1, 0)$ if $\rho< 0$; see also Example~\ref{exgaussian}
below. As a consequence, the cited results on the empirical copula
process do not apply to such copulas. This problem has been largely
ignored in the literature, and unjustified calls to the above results
abound. A~notable exception is the paper by Omelka, Gijbels, and Veraverbeke~\cite
{omelkagijbelsveraverbeke2009}. On page~3031 of that paper, it is
claimed that weak convergence of the empirical copula process still
holds if the first-order partial derivatives are continuous at $[0,
1]^2 \setminus\{(0, 0), (0, 1), (1, 0), (1, 1)\}$.

It is the aim of this paper to remedy the situation by showing that the
earlier cited results on the empirical copula process actually do hold
under a much less restrictive assumption, including indeed many copula
families that were hitherto excluded. The assumption is non-restrictive
in the sense that it is needed anyway to ensure that the candidate
limiting process exists and has continuous trajectories. The results
are stated and proved in general dimensions. When specialized to the
bivariate case, the condition is substantially weaker still than the
above-mentioned condition in Omelka, Gijbels, and Veraverbeke \cite{omelkagijbelsveraverbeke2009}.

Let $F$ be a $d$-variate cumulative distribution function (c.d.f.) with
continuous margins $F_1, \ldots, F_d$ and copula $C$, that is, $F(x) =
C(F_1(x_1), \ldots, F_d(x_d))$ for $x \in\reals^d$. Let $X_1, \ldots,
X_n$ be independent random vectors with common distribution $F$, where
$X_i = (X_{i1}, \ldots, X_{id})$. The empirical copula was defined in
Deheuvels \cite{deheuvels1979} as
%
\begin{equation}
\label{eempcop}
C_n(u) = F_n  ( F_{n1}^{-1}(u_1), \ldots, F_{nd}^{-1}(u_d)  ),\qquad
  u \in[0, 1]^d,
\end{equation}
where $F_n$ and $F_{nj}$ are the empirical joint and marginal cdfs of
the sample and where $F_{nj}^{-1}$ is the marginal quantile function of
the $j$th coordinate sample; see Section~\ref{secpreliminaries} below
for details. The empirical copula $C_n$ is invariant under monotone
increasing transformations on the data, so it depends on the data only
through the ranks. Indeed, up to a difference of order $1/n$, the
empirical copula can be seen as the empirical c.d.f. of the sample of
normalized ranks, as, for instance, in R{\"u}schendorf \cite{ruschendorf1976}. For
convenience, the definition in equation~\eqref{eempcop} will be
employed throughout the paper.

The empirical copula process is defined by
%
\begin{equation}
\label{eCCn}
\CC_n = \sqrt{n} (C_n - C),
\end{equation}
to be seen as a random function on $[0, 1]^d$. We are essentially
interested in the asymptotic distribution of $\CC_n$ in the space $\ell
^\infty([0, 1]^d)$ of bounded functions from $[0, 1]^d$ into $\reals$
equipped with the topology of uniform convergence. Weak convergence is
to be understood in the sense used in the monograph by van~der Vaart and Wellner \cite
{vandervaartwellner1996}, in particular their Definition~1.3.3.

Although the empirical copula is itself a rather crude estimator of
$C$, it plays a crucial rule in more sophisticated inference procedures
on $C$, much in the same way as the empirical c.d.f. $F_n$ is a
fundamental object for creating and understanding inference procedures
on $F$ or parameters thereof. For instance, the empirical copula is a
basic building block when estimating copula densities (Chen and Huang \cite
{chenhuang2007}, Omelka, Gijbels and Veraverbeke~\cite{omelkagijbelsveraverbeke2009}) or dependence
measures and functions (Schmid \textit{et~al.} \cite{schmidetal2010},
Genest and Segers \cite{genestsegers2010}), for testing for independence (Genest and R{\'e}millard~\cite
{genestremillard2004}, Genest, Quessy and R{\'e}millard~\cite{genestquessyremillard2007},
Kojadinovic and Holmes \cite{kojadinovicholmes2009}), for testing for shape constraints (Denuit and Scaillet \cite
{denuitscaillet2004}, Scaillet~\cite{scaillet2005}, Kojadinovic and Yan~\cite{kojadinovicyan2010}), for
resampling (R\'emillard and Scaillet \cite{remillardscaillet2009}, B{\"u}cher and Dette \cite{bucherdette2010}), and
so forth.

After some preliminaries in Section~\ref{secpreliminaries}, the
principal result of the paper is given in Section~\ref{secempproc},
stating weak convergence of the empirical copula process under the
condition that for every $j \in\{1, \ldots, d\}$, the $j$th
first-order partial derivative $\dot{C}_j$ exists and is continuous on
the set $\{ u \in[0, 1]^d \dvt 0 < u_j < 1 \}$. The condition is
non-restrictive in the sense that it is necessary for the candidate
limiting process to exist and have continuous trajectories. Moreover,
the resampling method based on the multiplier central limit theorem
proposed in R\'emillard and Scaillet \cite{remillardscaillet2009} is shown to be valid under
the same condition. Section~\ref{secstute} provides a refinement of
the main result: under certain bounds on the second-order partial
derivatives that allow for explosive behavior near the boundaries, the
almost sure error bound on the remainder term in Stute \cite{stute1984} and
Tsukahara \cite{tsukahara2005} can be entirely recovered. The result hinges on
an exponential inequality for a certain oscillation modulus of the
multivariate empirical process detailed in the \hyperref[app]{Appendix}; the inequality
is a generalization of a similar inequality in Einmahl \cite{einmahl1987} and
was communicated by Hideatsu Tsukahara. Section~\ref{secexamples}
concludes the paper with a number of examples of copulas that do or do
not verify certain sets of conditions.

\section{Preliminaries}
\label{secpreliminaries}

Let $X_i = (X_{i1}, \ldots, X_{id})$, $i \in\{1, 2, \ldots\}$, be
independent random vectors with common c.d.f.~$F$ whose margins $F_1,
\ldots, F_d$ are continuous and whose copula is denoted by~$C$. Define
$U_{ij} = F_j(X_{ij})$ for $i \in\{1, \ldots, n\}$ and $j \in\{1,
\ldots, d\}$. The random vectors $U_i = (U_{i1}, \ldots, U_{id})$
constitute an i.i.d. sample from $C$. Consider the following empirical
distribution functions: for $x \in\reals^d$ and for $u \in[0, 1]^d$,
\begin{eqnarray*}
F_n(x) &= &\frac{1}{n} \sum_{i=1}^n \ind_{(-\infty, x]}(X_i), \qquad
F_{nj}(x_j) = \frac{1}{n} \sum_{i=1}^n \ind_{(-\infty, x_j]}(X_{ij}),
\\
G_n(u) &=& \frac{1}{n} \sum_{i=1}^n \ind_{[0, u]}(U_i),\qquad
G_{nj}(u_j) = \frac{1}{n} \sum_{i=1}^n \ind_{[0, u_j]}(U_{ij}).
\end{eqnarray*}
Here, order relations on vectors are to be interpreted componentwise,
and $\ind_A(x)$ is equal to $1$ or~$0$ according to whether $x$ is an
element of $A$ or not. Let $X_{1:n,j} < \cdots< X_{n:n,j}$ and
$U_{1:n,j} < \cdots< U_{n:n,j}$ be the vectors of ascending order
statistics of the $j$th coordinate samples $X_{1j}, \ldots, X_{nj}$ and
$U_{1j}, \ldots, U_{nj}$, respectively. The marginal quantile functions
associated to $F_{nj}$ and $G_{nj}$ are
\begin{eqnarray*}
F_{nj}^{-1}(u_j)
&= &\inf\{ x \in\reals\dvt F_{nj}(x) \ge u_j \} \\
&=&
\cases{
X_{k:n,j}, &\quad $\mbox{if $(k-1)/n < u_j \le k/n$,}$\vspace*{2pt} \cr
-\infty,& \quad$\mbox{if $u_j = 0$;}$
}
\\
G_{nj}^{-1}(u_j)
&=& \inf\{ u \in[0, 1] \dvt G_{nj}(u) \ge u_j \} \\
&=&
\cases{
U_{k:n,j}, & \quad$\mbox{if $(k-1)/n < u_j \le k/n$,}$\vspace*{2pt}\cr
0, & \quad$\mbox{if $u_j = 0$.}$}
\end{eqnarray*}
Some thought shows that $X_{ij} \le F_{nj}^{-1}(u_j)$ if and only if
$U_{ij} \le G_{nj}^{-1}(u_j)$, for all $i \in\{1, \ldots, n\}$, $j \in
\{1, \ldots, d\}$ and $u_j \in[0, 1]$. It follows that the empirical
copula in equation~\eqref{eempcop} is given by
\[
C_n(u)
= G_n  ( G_{n1}^{-1}(u_1), \ldots, G_{nd}^{-1}(u_d)  ).
\]
In particular, without loss of generality we can work directly with the
sample $U_1, \ldots, U_n$ from~$C$.

The empirical processes associated to the empirical distribution
functions $G_n$ and $G_{nj}$ are given by
%
\begin{eqnarray}
\label{ealphan}
\alpha_n(u) = \sqrt{n}  \bigl( G_n(u) - C(u)  \bigr), \qquad
\alpha_{nj}(u_j) = \sqrt{n} \bigl( G_{nj}(u_j) - u_j  \bigr),
\end{eqnarray}
for $u \in[0, 1]^d$ and $u_j \in[0, 1]$. Note that $\alpha_{nj}(0) =
\alpha_{nj}(1) = 0$ almost surely. We have
\[
\alpha_n \dto\alpha\qquad (n \to\infty)
\]
in $\ell^\infty([0, 1]^d)$, the arrow `$\dto$' denoting weak
convergence as in Definition~1.3.3 in van~der Vaart and Wellner \cite{vandervaartwellner1996}.
The limit process $\alpha$ is a $C$-Brownian bridge, that is, a tight
Gaussian process, centered and with covariance function
\[
\cov ( \alpha(u), \alpha(v)  ) = C ( u \wedge v ) - C(u)   C(v),
\]
for $u, v \in[0, 1]^d$; here $u \wedge v = (\min(u_1, v_1), \ldots,
\min(u_d, v_d))$. Tightness of the process~$\alpha$ and continuity of
its mean and covariance functions implies the existence of a version of~$\alpha$ with continuous trajectories. Without loss of generality, we
assume henceforth that~$\alpha$ is such a version.

For $j \in\{1, \ldots, d\}$, let $e_j$ be the $j$th coordinate vector
in $\reals^d$. For $u \in[0, 1]^d$ such that $0 < u_j < 1$, let
\[
\dot{C}_j(u) = \lim_{h \to0} \frac{C(u + he_j) - C(u)}{h},
\]
be the $j$th first-order partial derivative of $C$, provided it exists.

\begin{condition}
\label{cdiffC}
For each $j \in\{1, \ldots, d\}$, the $j$th first-order partial
derivative $\dot{C}_j$ exists and is continuous on the set $V_{d,j} :=
\{ u \in[0, 1]^d \dvt 0 < u_j < 1 \}$.
\end{condition}

Henceforth, assume Condition~\ref{cdiffC} holds. To facilitate
notation, we will extend the domain of $\dot{C}_j$ to the whole of $[0,
1]^d$ by setting
%
\begin{equation}
\label{eextend}
\dot{C}_j(u) =
\cases{
 \displaystyle\limsup_{h \downarrow0} \frac{C(u + he_j)}{h}, & \quad$\mbox
{if $u \in[0, 1]^d$, $u_j = 0$,}$\vspace*{2pt}\cr
 \displaystyle\limsup_{h \downarrow0} \frac{C(u) - C(u - he_j)}{h}, &\quad
$\mbox{if $u \in[0, 1]^d$, $u_j = 1$.}$}
\end{equation}
In this way, $\dot{C}_j$ is defined everywhere on $[0, 1]^d$, takes
values in $[0, 1]$ (because $|C(u) - C(v)| \le\sum_{j=1}^d |u_j -
v_j|$), and is continuous on the set $V_{d,j}$, by virtue of
Condition~\ref{cdiffC}. Also note that $\dot{C}_j(u) = 0$ as soon as
$u_i = 0$ for some $i \ne j$. 

\section{Weak convergence}
\label{secempproc}

In Proposition~\ref{pempproc}, Condition~\ref{cdiffC} is shown to be
sufficient for the weak convergence of the empirical copula process $\CC
_n$. In contrast to earlier results, Condition~\ref{cdiffC} does not
require existence or continuity of the partial derivatives on certain
boundaries. Although the improvement is seemingly small, it
dramatically enlarges the set of copulas to which it applies; see
Section~\ref{secexamples}. Similarly, the unconditional multiplier
central limit theorem for the empirical copula process based on
estimated first-order partial derivatives continues to hold
(Proposition~\ref{pmclt}). This result is useful as a justification of
certain resampling procedures that serve to compute critical values for
test statistics based on the empirical copula in case of a composite
null hypothesis, for instance, in the context of goodness-of-fit
testing as in~Kojadinovic and Yan \cite{kojadinovicyan2010}.

Assume first that the first-order partial derivatives $\dot{C}_j$ exist
and are continuous throughout the closed hypercube $[0, 1]^d$. For $u
\in[0, 1]^d$, define
%
\begin{equation}
\label{eCC}
\CC(u) = \alpha(u) - \sum_{j=1}^d \dot{C}_j(u)   \alpha_j(u_j),
\end{equation}
where $\alpha_j(u_j) = \alpha(1, \ldots, 1, u_j, 1, \ldots, 1)$, the
variable $u_j$ appearing at the $j$th entry. By continuity of $\dot
{C}_j$ throughout $[0, 1]^d$, the trajectories of $\CC$ are continuous.
From Fermanian \textit{et al.} \cite{fermanianradulovicwegkamp2004} and Tsukahara \cite{tsukahara2005}
we learn that $\CC_n \dto\CC$ as $n \to\infty$ in the space $\ell
^\infty([0, 1]^d)$.

The structure of the limit process $\CC$ in equation~\eqref{eCC} can
be understood as follows. The first term, $\alpha(u)$, would be there
even if the true margins $F_j$ were used rather than their empirical
counterparts $F_{nj}$. The terms $- \dot{C}_j(u)   \alpha_j(u_j)$
encode the impact of not knowing the true quantiles $F_j^{-1}(u_j)$ and
having to replace them by the empirical quantiles $F_{nj}^{-1}(u_j)$.
The minus sign comes from the Bahadur--Kiefer result stating that $\sqrt
{n}(G_{nj}^{-1}(u_j) - u_j)$ is asymptotically undistinguishable from
$- \sqrt{n} (G_{nj}(u_j) - u_j)$; see, for instance, Shorack and Wellner \cite
{shorackwellner1986}, Chapter~15. The partial derivative $\dot
{C}_j(u)$ quantifies the sensitivity of $C$ with respect to small
deviations in the $j$th margin.

Now consider the same process $\CC$ as in equation~\eqref{eCC} but
under Condition~\ref{cdiffC} and with the domain of the partial
derivatives extended to $[0, 1]^d$ as in equation~\eqref{eextend}.
Since the trajectories of~$\alpha$ are continuous and since $\alpha
_j(0) = \alpha_j(1) = 0$ for each $j \in\{1, \ldots, d\}$, the
trajectories of $\CC$ are continuous, even though $\dot{C}_j$ may fail
to be continuous at points $u \in[0, 1]^d$, such that $u_j \in\{0, 1\}
$. The process $\CC$ is the weak limit in $\ell^\infty([0, 1]^d)$ of
the sequence of processes
%
\begin{equation}
\label{eCCntilde}
\tilde{\CC}_n(u) = \alpha_n(u) - \sum_{j=1}^d \dot{C}_j(u)   \alpha
_{nj}(u_j),  \qquad u \in[0, 1]^d.
\end{equation}
The reason is that the map from $\ell^\infty([0, 1]^d)$ into itself
that sends a function $f$ to $f - \sum_{j = 1}^d \dot{C}_j   \pi
_j(f)$, where $(\pi_j(f))(u) = f(1, \ldots, 1, u_j, 1, \ldots, 1)$, is
linear and bounded.

\begin{proposition}
\label{pempproc}
If Condition~\ref{cdiffC} holds, then, with $\tilde{\CC}_n$ as in
equation~\eqref{eCCntilde},
\[
\sup_{u \in[0, 1]^d}  | \CC_n(u) - \tilde{\CC}_n(u)  | \pto0\qquad
 (n \to\infty).
\]
As a consequence, in $\ell^\infty([0, 1]^d)$,
\[
\CC_n \dto\CC\qquad (n \to\infty).
\]
\end{proposition}

\begin{pf}
It suffices to show the first statement of the proposition. For $u \in
[0, 1]^d$, put
\[
R_n(u) =  | \CC_n(u) - \tilde{\CC}_n(u)  |,  \qquad u \in[0, 1]^d.
\]
If $u_j = 0$ for some $j \in\{1, \ldots, d\}$, then obviously $\CC
_n(u) = \tilde{\CC}_n(u) = 0$, so $R_n(u) = 0$ as well. The vector of
marginal empirical quantiles is denoted by
%
\begin{equation}
\label{evn}
v_n(u) =  ( G_{n1}^{-1}(u_1), \ldots, G_{nd}^{-1}(u_d)  ),
 \qquad u \in[0, 1]^d.
\end{equation}
We have
%
\begin{eqnarray}
\label{edecomp1}
\CC_n(u) &= &\sqrt{n}  \bigl( C_n(u) - C(u)  \bigr) \nonumber\\
&= &\sqrt{n}
 \{
G_n  ( v_n(u)  ) - C  ( v_n(u)  )
 \}
+ \sqrt{n}
 \{
C  ( v_n(u)  ) - C(u)
 \} \\
&=& \alpha_n  ( v_n(u)  )
+ \sqrt{n}
 \{
C  ( v_n(u)  ) - C(u)
 \}.\nonumber
\end{eqnarray}
Since $\alpha_n$ converges weakly in $\ell^\infty([0, 1]^d)$ to a
$C$-Brownian bridge $\alpha$, whose trajectories are continuous, the
sequence $(\alpha_n)_n$ is asymptotically uniformly equicontinuous; see
Theorem~1.5.7 and Addendum~1.5.8 in van~der Vaart and Wellner \cite{vandervaartwellner1996}. As
$\sup_{u_j \in[0, 1]} | G_{nj}^{-1}(u_j) - u_j | \to0$ almost surely,
it follows that
\[
\sup_{u \in[0, 1]^d}
 |
\alpha_n  ( v_n(u)  )
- \alpha_n(u)
 |
\pto0 \qquad (n \to\infty).
\]

Fix $u \in[0, 1]^d$. Put $w(t) = u + t \{ v_n(u) - u \}$ and $f(t) =
C(w(t))$ for $t \in[0, 1]$. If $u \in(0, 1]^d$, then $v_n(u) \in(0,
1)^d$, and therefore $w(t) \in(0, 1)^d$ for all $t \in(0, 1]$, as
well. By Condition~\ref{cdiffC}, the function $f$ is continuous on
$[0, 1]$ and continuously differentiable on $(0, 1)$. By the mean value
theorem, there exists $t^* = t_n(u) \in(0, 1)$ such that $f(1) - f(0)
= f'(t^*)$, yielding
%
\begin{equation}
\label{edecomp2}
\sqrt{n}  \{ C  ( v_n(u)  ) - C(u)  \}
= \sum_{j=1}^d \dot{C}_j ( w(t^*)  )   \sqrt{n}    \bigl(
G_{nj}^{-1}(u_j) - u_j  \bigr).
\end{equation}
If one or more of the components of $u$ are zero, then the above
display remains true as well, no matter how $t^* \in(0, 1)$ is
defined, because both sides of the equation are equal to zero. In
particular, if $u_k = 0$ for some $k \in\{1, \ldots, d\}$, then the
$k$th term on the right-hand side vanishes because $G_{nk}^{-1}(0) = 0$
whereas the terms with index $j \ne k$ vanish because the $k$th
component of the vector $w(t^*)$ is zero, and thus the first-order
partial derivatives $\dot{C}_j$ vanish at this point.

It is known since Kiefer \cite{kiefer1970} that
\[
\sup_{u_j \in[0, 1]}  \bigl| \sqrt{n}  \bigl( G_{nj}^{-1}(u_j) - u_j
 \bigr) + \alpha_{nj}(u_j)  \bigr| \pto0\qquad  (n \to\infty).
\]
Since $0 \le\dot{C}_j \le1$, we find
\[
\sup_{u \in[0, 1]^d}
\Biggl|
\sqrt{n}  \{ C  ( v_n(u)  ) - C(u)  \}
+ \sum_{j=1}^d \dot{C}_j \bigl( u + t^* \{v_n(u) - u\}  \bigr) \alpha_{nj}(u_j)
\Biggr|
\pto0
\]
as $n \to\infty$. It remains to be shown that
\[
\sup_{u \in[0, 1]^d} D_{nj}(u) \pto0\qquad  (n \to\infty)
\]
for all $j \in\{1, \ldots, d\}$, where
%
\begin{equation}
\label{eDnj}
D_{nj}(u) =  \bigl| \dot{C}_j \bigl( u + t^* \{v_n(u) - u\}  \bigr) - \dot
{C}_j(u)  \bigr|    | \alpha_{nj}(u_j)  |.
\end{equation}
Fix $\eps> 0$ and $\delta\in(0, 1/2)$. Split the supremum over $u
\in[0, 1]^d$ according to the cases $u_j \in[\delta, 1 - \delta]$ on
the one hand and $u_j \in[0, \delta) \cup(1-\delta, 1]$ on the other
hand. We have
\begin{eqnarray*}
\Pr\Bigl( \sup_{u \in[0, 1]^d} D_{nj}(u) > \eps\Bigr)
&\le&\Pr\Bigl( \sup_{u \in[0, 1]^d, u_j \in[\delta, 1-\delta]}
D_{nj}(u) > \eps\Bigr) \\
&&{}+ \Pr\Bigl( \sup_{u \in[0, 1]^d, u_j \notin[\delta, 1-\delta]}
D_{nj}(u) > \eps\Bigr).
\end{eqnarray*}
Since $\sup_{u \in[0, 1]^d} |v_n(u) - u| \to0$ almost surely, since
$\dot{C}_j$ is uniformly continuous on $\{ u \in[0, 1]^d\dvt \delta/2 \le
u_j \le1 - \delta/2 \}$, and since the sequence $\sup_{u_{nj} \in[0,
1]} |\alpha_{nj}(u_j)|$ is bounded in probability, the first
probability on the right-hand side of the previous display converges to
zero. As $|x - y| \le1$ whenever $x, y \in[0, 1]$ and since $0 \le
\dot{C}_j(w) \le1$ for all $w \in[0, 1]^d$, the second probability on
the right-hand side of the previous display is bounded by\looseness=1
\[
\Pr\Bigl( \sup_{u_j \in[0, \delta) \cup(1-\delta, 1]}  | \alpha
_{nj}(u_j)  | > \eps\Bigr).
\]\looseness=0
By the portmanteau lemma, the $\limsup$ of this probability as $n \to
\infty$ is bounded by
\[
\Pr\Bigl( \sup_{u_j \in[0, \delta) \cup(1-\delta, 1]}  | \alpha
_{j}(u_j)  | \ge\eps\Bigr).
\]
The process $\alpha_j$ being a standard Brownian bridge, the above
probability can be made smaller than an arbitrarily chosen $\eta> 0$
by choosing $\delta$ sufficiently small. We find
\[
\limsup_{n \to\infty} \Pr\Bigl( \sup_{u \in[0, 1]^d} D_{nj}(u) >
\eps\Bigr) \le\eta.
\]
As $\eta$ was arbitrary, the claim is proven.\vspace*{3pt}
\end{pf}

An alternative to the direct proof above is to invoke the functional
delta method as in Fermanian \textit{et al.} \cite{fermanianradulovicwegkamp2004}. Required
then is a generalization of Lemma~2 in the cited paper asserting
Hadamard differentiability of a certain functional under Condition~\ref{cdiffC}.
This program is carried out for the bivariate case in B\"{u}cher \cite{buecher2011},
Lemma~2.6.

For purposes of hypothesis testing or confidence interval construction,
resampling procedures are often required; see the references in the
introduction. In Fermanian \textit{et al.} \cite{fermanianradulovicwegkamp2004}, a bootstrap
procedure for the empirical copula process is proposed, whereas in R\'emillard and Scaillet \cite
{remillardscaillet2009}, a method based on the multiplier central
limit theorem is employed. Yet another method is proposed in B{\"u}cher and Dette \cite
{bucherdette2010}. In the latter paper, the finite-sample properties
of all these methods are compared in a simulation study, and the
multiplier approach by R\'emillard and Scaillet \cite{remillardscaillet2009} is found to be
best overall. Although the latter approach requires estimation of the
first-order partial derivatives, it remains valid under Condition~\ref
{cdiffC}, allowing for discontinuities on the boundaries.\looseness=1

Let $\xi_1, \xi_2, \ldots$ be an i.i.d. sequence of random variables,
independent of the random vectors $X_1, X_2, \ldots,$ and with zero
mean, unit variance, and such that $\int_0^\infty\sqrt{\Pr(|\xi_1| >
x)}   \,\diff x < \infty$. Define
%
\begin{equation}
\label{ealphanprime}
\alpha_n'(u)
= \frac{1}{\sqrt{n}}
\sum_{i=1}^n \xi_i
 \bigl(
\ind\{ X_{i1} \le F_{n1}^{-1}(u_1), \ldots, X_{id} \le
F_{nd}^{-1}(u_d) \} - C_n(u)
 \bigr).
\end{equation}
In $ ( \ell^{\infty}([0, 1]^d)  )^2$, we have by Lemma~A.1 in
R\'emillard and Scaillet \cite{remillardscaillet2009},
%
\begin{equation}
\label{emcltalpha}
(\alpha_n, \alpha_n') \dto(\alpha, \alpha')\qquad  (n \to\infty),
\end{equation}
where $\alpha'$ is an independent copy of $\alpha$. Further, let $\hat
{\dot{C}}_{nj}(u)$ be an estimator of $\dot{C}_j(u)$; for instance,
apply finite differencing to the empirical copula at a spacing
proportional to $n^{-1/2}$ as in R\'emillard and Scaillet \cite{remillardscaillet2009}. Define
%
\begin{equation}
\label{eCCnprime}
\CC_n'(u)
= \alpha_n'(u) - \sum_{j=1}^d \hat{\dot{C}}_{nj}(u)   \alpha_{nj}'(u_j),
\end{equation}
where $\alpha_{nj}'(u_j) = \alpha_n'(1, \ldots, 1, u_j, 1, \ldots, 1)$,
the variable $u_j$ appearing at the $j$th coordinate.

\begin{proposition}
\label{pmclt}
Assume Condition~\ref{cdiffC}. If there exists a constant $K$ such
that $|\hat{\dot{C}}_{nj}(u)| \le K$ for all $n, j, u$, and if
%
\begin{equation}
\label{ediffCestim}
\sup_{u \in[0, 1]^d : u_j \in[\delta, 1-\delta]}  | \hat{\dot
{C}}_{nj}(u) - \dot{C}_j(u)  | \pto0
\qquad (n \to\infty)
\end{equation}
for all $\delta\in(0, 1/2)$ and all $j \in\{1, \ldots, d\}$, then in
$ ( \ell^{\infty}([0, 1]^d)  )^2$, we have
\[
(\CC_n, \CC_n') \dto(\CC, \CC')\qquad  (n \to\infty),
\]
where $\CC'$ is an independent copy of $\CC$.
\end{proposition}

\begin{pf}
Recall the process $\alpha_n'$ in equation~\eqref{ealphanprime}, and define
\[
\tilde{\CC}_n'(u) = \alpha_n'(u) - \sum_{j=1}^d \dot{C}_j(u)   \alpha
_{nj}'(u_j),  \qquad u \in[0, 1]^d.
\]
The difference with the process $\CC_n'$ in equation~\eqref{eCCnprime}
is that the true partial derivatives of $C$ are used rather than the
estimated ones. By Proposition~\ref{pempproc} and equation~\eqref
{emcltalpha}, we have
\[
(\CC_n, \tilde{\CC}_n') \dto(\CC, \CC')\qquad  (n \to\infty)
\]
in $ ( \ell^{\infty}([0, 1]^d)  )^2$. Moreover,
\[
 | \CC_n'(u) - \tilde{\CC}_n'(u)  |
\le\sum_{j=1}^d  | \hat{\dot{C}}_{nj}(u) - \dot{C}_j(u)  |
| \alpha_{nj}'(u_j) |.
\]
It suffices to show that each of the $d$ terms on the right-hand side
converges to $0$ in probability, uniformly in $u \in[0, 1]^d$. The
argument is similar to the one at the end of the proof of
Proposition~\ref{pempproc}. Pick $\delta\in(0, 1/2)$, and split the
supremum according to the cases $u_j \in[\delta, 1-\delta]$ and $u_j
\in[0, \delta) \cup(1-\delta, 1]$. For the first case, use
equation~\eqref{ediffCestim} together with tightness of $\alpha
_{nj}'$. For the second case, use the assumed uniform boundedness of
the partial derivative estimators and the fact that the limit process
$\hat{\alpha}_j$ is a standard Brownian bridge, having continuous
trajectories and vanishing at $0$ and $1$.
\end{pf}

\section{Almost sure rate}
\label{secstute}

Recall the empirical copula process $\CC_n$ in equation~\eqref{eCCn}
together with its approximation~$\tilde{\CC}_n$ in equation~\eqref
{eCCntilde}. If the second-order partial derivatives of $C$ exist and
are continuous on $[0, 1]^d$, then the original result by Stute \cite
{stute1984}, proved in detail in Tsukahara~\cite{tsukahara2000}, reinforces the
first claim of Proposition~\ref{pempproc} to
%
\begin{eqnarray}
\label{estute}
&&\sup_{u \in[0, 1]^d} |\CC_n(u) - \tilde{\CC}_n(u)|
\nonumber
\\[-8pt]
\\[-8pt]
\nonumber
&&\quad= \mathrm{O}  ( n^{-1/4} (\log n)^{1/2} (\log\log n)^{1/4}  )\qquad
(n \to\infty)  \mbox{ almost surely.}
\end{eqnarray}
For many copulas, however, the second-order partial derivatives explode
near certain parts of the boundaries. The question then is how this
affects the above rate. Recall $V_{d,j} = \{ u \in[0, 1]^d \dvt 0 < u_j <
1 \}$ for $j \in\{1, \ldots, d\}$.

\begin{condition}
\label{cdiffCK}
For every $i, j \in\{1, \ldots, d\}$, the second-order partial
derivative $\ddot{C}_{ij}$ is defined and continuous on the set
$V_{d,i} \cap V_{d,j}$, and there exists a constant $K > 0$ such that
\[
|\ddot{C}_{ij}(u)| \le K   \min\biggl( \frac{1}{u_i(1-u_i)} ,   \frac
{1}{u_j(1-u_j)} \biggr),  \qquad u \in V_{d,i} \cap V_{d,j}.
\]
\end{condition}

Condition~\ref{cdiffCK} holds, for instance, for absolutely continuous
bivariate Gaussian copulas and for bivariate extreme-value copulas
whose Pickands dependence functions are twice continuously
differentiable and satisfy a certain bound; see Section~\ref{secexamples}.

Under Condition~\ref{cdiffCK}, the rate in equation~\eqref{estute}
can be entirely recovered. The following proposition has benefited from
a suggestion of John H.J. Einmahl leading to an improvement of a
result in an earlier version of the paper claiming a slightly slower
rate. Furthermore, part of the proof is an adaptation due to Hideatsu
Tsukahara of the end of the proof of Theorem~4.1 in Tsukahara \cite
{tsukahara2000}, upon which the present result is based.

\begin{proposition}
\label{pempprocrate}
If Conditions~\ref{cdiffC} and~\ref{cdiffCK} are verified, then
equation~\eqref{estute} holds.
\end{proposition}

\begin{pf}
Combining equations~\eqref{edecomp1} and \eqref{edecomp2} in the
proof of Proposition~\ref{pempproc} yields
\[
\CC_n(u) = \alpha_n (v_n(u) ) + \sum_{j=1}^d \dot{C}_j (
w(t^*)  )   \sqrt{n}    \bigl( G_{nj}^{-1}(u_j) - u_j  \bigr),
\qquad  u \in[0, 1]^d,
\]
with $\alpha_n$ the ordinary multivariate empirical process in
equation~\eqref{ealphan}, $v_n(u)$ the vector of marginal empirical
quantiles in equation~\eqref{evn}, and $w(t^*) = u + t^* \{ v_n(u) - u
\}$ a certain point on the line segment between $u$ and $v_n(u)$ with
local coordinate $t^* = t_n(u) \in(0, 1)$. In view of the definition
of $\tilde{\CC}_n(u)$ in equation~\eqref{eCCntilde}, it follows that
\[
\sup_{u \in[0, 1]^d} |\CC_n(u) - \tilde{\CC}_n(u)| \le\I_n + \II_n +
\III_n,
\]
where
\begin{eqnarray*}
\I_n &=& \sup_{u \in[0, 1]^d}  | \alpha_n  ( v_n(u)  ) -
\alpha_n(u)  |, \\
\II_n &= &\sum_{j=1}^d \sup_{u \in[0, 1]^d}  \bigl| \sqrt{n}  \bigl(
G_{nj}^{-1}(u_j) - u_j  \bigr) + \alpha_{nj}(u_j)  \bigr|, \\
\III_n &=& \sum_{j=1}^d \sup_{u \in[0, 1]^d} D_{nj}(u),
\end{eqnarray*}
with $D_{nj}(u)$ as defined in equation~\eqref{eDnj}. By Kiefer \cite
{kiefer1970}, the term $\II_n$ is $\mathrm{O}  ( n^{-1/4} (\log n)^{1/2}
\times (\log\log n)^{1/4}  )$ as $n \to\infty$, almost surely. It
suffices to show that the same almost sure rate is valid for $\I_n$ and
$\III_n$, too.

\textit{The term $\I_n$.}
The argument is adapted from the final part of the proof of Theorem~4.1
in Tsukahara \cite{tsukahara2000}, and its essence was kindly provided by
Hideatsu Tsukahara. We have
\[
\I_n \le M_n (A_n), \qquad A_n = \max_{j \in\{1, \ldots, d\}} \sup_{u_j
\in[0, 1]} | G_{nj}^{-1}(u_j) - u_j |,
\]
and $M_n(a)$ is the oscillation modulus of the multivariate empirical
process $\alpha_n$ defined in equation~\eqref{eoscillation}. We will
employ the exponential inequality for $\Pr\{ M_n(a) \ge\lambda\}$
stated in Proposition~\ref{poscillation}, which generalizes
Inequality~3.5 in Einmahl \cite{einmahl1987}. Set $a_n = n^{-1/2} (\log\log
n)^{1/2}$. By the Chung--Smirnov law of the iterated logarithm for
empirical distribution functions (see, e.g., Shorack and Wellner \cite
{shorackwellner1986}, page 504),
%
\begin{eqnarray}
\label{eLIL}
\limsup_{n \to\infty} \frac{1}{a_n} \sup_{u_j \in[0, 1]} |
G_{nj}^{-1}(u_j) - u_j |
&= &\limsup_{n \to\infty} \frac{1}{a_n} \sup_{v_j \in[0, 1]} | v_j -
G_{nj}(v_j) |
\nonumber
\\[-8pt]
\\[-8pt]
\nonumber
&= &1/\sqrt{2}\qquad  \mbox{almost surely}.
\end{eqnarray}
Choose $\lambda_n = 2 K_2^{-1/2}   n^{-1/4} (\log n)^{1/2} (\log\log
n)^{1/4}$ for $K_2$ as in Proposition~\ref{poscillation}. Since
$\lambda_n / (n^{1/2} a_n) \to0$ as $n \to\infty$, and since the
function $\psi$ in equation~\eqref{epsi} below is decreasing with $\psi
(0) = 1$, it follows that $\psi(\lambda_n / (n^{1/2} a_n)) \ge1/2$ for
sufficiently large $n$. Furthermore, we have
\[
\sum_{n \ge2} \frac{1}{a_n}   \exp\biggl( - \frac{K_2 \lambda
_n^2}{2a_n} \biggr) = \sum_{n \ge2} \frac{1}{n^{3/2} (\log\log
n)^{1/2}} < \infty.
\]
By the Borel--Cantelli lemma and Proposition~\eqref{poscillation}, as
$n \to\infty$,
\[
\I_n \le M_n(A_n) \le M_n(a_n) = \mathrm{O}  ( n^{-1/4} (\log n)^{1/2} (\log
\log n)^{1/4}  )\qquad  \mbox{almost surely}.
\]

\textit{The term $\III_n$.}
Let
\[
\delta_n = n^{-1/2}   (\log n)   (\log\log n)^{-1/2}.
\]
Fix $j \in\{1, \ldots, d\}$. We split the supremum of $D_{nj}(u)$ over
$u \in[0, 1]^d$ according to the cases $u_j \in[0, \delta_n) \cup
(1-\delta_n, 1]$ and $u_j \in[\delta_n, 1-\delta_n]$.

Since $0 \le\dot{C}_j \le1$, the supremum over $u \in[0, 1]^d$ such
that $u_j \in[0, \delta_n) \cup(1-\delta_n, 1]$ is bounded by
\[
\sup_{u \in[0, 1]^d: u_j \in[0, \delta_n) \cup(1-\delta_n, 1]}
D_{nj}(u) \le\sup_{u_j \in[0, \delta_n) \cup(1-\delta_n, 1]} |\alpha
_{nj}(u_j)|.
\]
By Theorem~2.(iii) in Einmahl and Mason \cite{einmahlmason1988} applied to $(d, \nu,
k_n) = (1, 1/2, n   \delta_n)$, the previous supremum is of the order
%
\begin{eqnarray}
\label{ebound1}
\sup_{u_j \in[0, \delta_n) \cup(1-\delta_n]} |\alpha_{nj}(u_j)| &= &\mathrm{O}
 ( \delta_n^{1/2} (\log\log n)^{1/2}  )
\nonumber
\\[-8pt]
\\[-8pt]
\nonumber
&=&\mathrm{O}  ( n^{-1/4}   (\log n)^{1/2}  (\log\log n)^{1/4}  ) \qquad
(n \to\infty)  \mbox{ almost surely.}\qquad
\end{eqnarray}

Next let $u \in[0, 1]^d$ be such that $\delta_n \le u_j \le1-\delta
_n$. By Lemma~\ref{lincr} below and by convexity of the function $(0,
1) \ni s \mapsto1/\{s(1-s)\}$,
\begin{eqnarray*}
D_{nj}(u)
&= & \bigl| \dot{C}_j \bigl( u + \lambda_n(u) \{v_n(u) - u\}  \bigr) - \dot
{C}_j(u)  \bigr|    | \alpha_{nj}(u_j)  | \\
&\le& K   \max\biggl( \frac{1}{u_j(1-u_j)} ,   \frac
{1}{G_{nj}^{-1}(u_j) (1 - G_{nj}^{-1}(u_j))} \biggr)   \norm{ v_n(u) -
u }_1    | \alpha_{nj}(u_j)  |,
\end{eqnarray*}
with $\norm{x}_1 = \sum_{j=1}^d |x_j|$. Let $b_n = (\log n)^{1/2}
\log\log n$; clearly $\sum_{n=2}^\infty n^{-1} b_n^{-2} < \infty$. By
Cs{\'a}ki~\cite{csaki1974} or Mason \cite{mason1981},
\[
\Pr\biggl( \sup_{0 < s < 1} \frac{|\alpha_{nj}(s)|}{(s(1-s))^{1/2}} >
b_n \mbox{ infinitely often} \biggr) = 0.
\]
It follows that, with probability one, for all sufficiently large $n$,
\[
|\alpha_{nj}(u_j)| \le \bigl( u_j(1-u_j)  \bigr)^{1/2}   b_n,\qquad
u_j \in[0, 1].
\]
Let $I$ denote the identity function on $[0, 1]$, and let $\norm{ \cdot
 }_\infty$ denote the supremum norm. For $u_j \in[\delta_n, 1-\delta_n]$,
\begin{eqnarray*}
G_{nj}^{-1}(u_j) &=& u_j  \biggl ( 1 + \frac{G_{nj}^{-1}(u_j) -
u_j}{u_j} \biggr)\ge u_j  \biggl ( 1 - \frac{\norm{G_{nj}^{-1} - I}_\infty}{\delta_n}
\biggr), \\
1 - G_{nj}^{-1}(u_j) &\ge&(1-u_j) \biggl( 1 - \frac{\norm{G_{nj}^{-1}
- I}_\infty}{\delta_n} \biggr).
\end{eqnarray*}
By the law of the iterated logarithm (see \eqref{eLIL})
\[
\norm{G_{nj}^{-1} - I} = \mathrm{o}(\delta_n)\qquad  (n \to\infty)  \mbox{ almost surely.}
\]
We find that with probability one, for all sufficiently large $n$ and
for all $u \in[0, 1]^d$ such that $u_j \in[\delta_n, 1 - \delta_n]$,
\[
D_{nj}(u) \le2K    \bigl(u_j(1-u_j) \bigr)^{-1/2}   \norm{v_n(u) -
u}_1   b_n.
\]
We use again the law of the iterated logarithm in \eqref{eLIL} to
bound $\norm{ v_n(u) - u }_1$. As a~consequence, with probability one,
%
\begin{eqnarray}
\label{ebound2}
\sup_{u \in[0, 1]^d : u_j \in[\delta_n, 1-\delta_n]} D_{nj}(u)
&=& \mathrm{O} ( \delta_n^{-1/2}   (\log\log n)^{1/2}   n^{-1/2}   b_n
 )
 \nonumber
 \\[-8pt]
 \\[-8pt]
 \nonumber
&=& \mathrm{O}  ( n^{-1/4}   (\log\log n)^{7/4}  )\qquad
 (n \to\infty)  \mbox{ almost surely}.
\end{eqnarray}
The bound in \eqref{ebound2} is dominated by the one in \eqref
{ebound1}. The latter therefore yields the total rate.
\end{pf}

\begin{lemma}
\label{lincr}
If Conditions~\ref{cdiffC} and~\ref{cdiffCK} hold, then
%
\begin{equation}
\label{eincr}
|\dot{C}_j(v) - \dot{C}_j(u)| \le K   \max\biggl( \frac
{1}{u_j(1-u_j)} ,   \frac{1}{v_j(1-v_j)} \biggr)   \norm{v - u}_1,
\end{equation}
for every $j \in\{1, \ldots, d\}$ and for every $u, v \in[0, 1]^d$
such that $0 < u_j < 1$ and $0 < v_j < 1$; here $\norm{x}_1 = \sum
_{i=1}^d |x_i|$ denotes the $L_1$-norm.
\end{lemma}

\begin{pf}
Fix $j \in\{1, \ldots, d\}$ and $u, v \in[0, 1]^d$ such that $u_j,
v_j \in(0, 1)$. Consider the line segment $w(t) = u + t(v-u)$ for $t
\in[0, 1]$, connecting $w(0) = u$ with $w(1) = v$; put $w_i(t) = u_i +
t(v_i - u_i)$ for $i \in\{1, \ldots, d\}$. Clearly $0 < w_j(t) < 1$
for all $t \in[0, 1]$. Next, consider the function $f(t) = \dot
{C}_j(w(t))$ for $t \in[0, 1]$. The function $f$ is continuous on $[0,
1]$ and continuously differentiable on $(0, 1)$. Indeed, if $u_i \ne
v_i$ for some $i \in\{1, \ldots, d\}$, then $0 < w_i(t) < 1$ for all
$t \in(0, 1)$; if $u_i = v_i$, then $w_i(t) = u_i = v_i$ does not
depend on $t$ at all. Hence, the derivative of $f$ in $t \in(0, 1)$ is
given by
\[
f'(t) = \sum_{i \in I} (v_i - u_i)   \ddot{C}_{ij}(w(t)),
\]
where $\mathcal{I} = \{ i \in\{1, \ldots, d\} \dvt u_i \ne v_i \}$. By
the mean-value theorem, we obtain that for some $t^* \in(0, 1)$,
\[
\dot{C}_j(v) - \dot{C}_j(u) = f(1) - f(0) = f'(t^*) = \sum_{i \in
\mathcal{I}} (v_i - u_i)   \ddot{C}_{ij}(w(t^*)).
\]
As a consequence,
\[
|\dot{C}_j(u) - \dot{C}_j(v)| \le\norm{v-u}_1   \max_{i \in\mathcal
{I}} \sup_{0 < t < 1} |\ddot{C}_{ij}(w(t))|.
\]
By Condition~\ref{cdiffCK},
\[
|\dot{C}_j(u) - \dot{C}_j(v)| \le\norm{v-u}_1   K   \sup_{0 < t < 1}
\frac{1}{w_j(t)   \{ 1 - w_j(t) \}}.
\]
Finally, since the function $s \mapsto1 / \{ s (1-s) \}$ is convex on
$(0, 1)$ and since $w_j(t)$ is a convex combination of $u_j$ and $v_j$,
the supremum of $1 / [ w_j(t) \{ 1 - w_j(t) \} ]$ over $t \in[0, 1]$
must be attained at one of the endpoints $u_j$ or $v_j$. Equation~\eqref
{eincr} follows.
\end{pf}

\section{Examples}
\label{secexamples}

\begin{example}[(Gaussian copula)]
\label{exgaussian}
Let $C$ be the $d$-variate Gaussian copula with correlation matrix $R
\in\reals^{d \times d}$, that is,
\[
C(u) = \Pr\Biggl( \bigcap_{j=1}^d \{ \Phi(X_j) \le u_j \} \Biggr),
\qquad  u \in[0, 1]^d,
\]
where $X = (X_1, \ldots, X_d)$ follows a $d$-variate Gaussian
distribution with zero means, unit variances, and correlation matrix
$R$; here $\Phi$ is the standard normal c.d.f. It can be checked readily
that if the correlation matrix $R$ is of full rank, then Condition~\ref
{cdiffC} is verified, and Propositions~\ref{pempproc} and~\ref
{pmclt} apply.

Still, if $0 < \rho_{1j} = \corr(X_1, X_j) < 1$ for all $j \in\{2,
\ldots, d\}$, then on the one hand we have $\lim_{u_1 \downarrow0} \dot
{C}_1(u_1, u_{-1}) = 1$ for all $u_{-1} \in(0, 1]^{d-1}$, whereas on
the other hand we have $\dot{C}_1(u) = 0$ as soon as $u_j = 0$ for some
$j \in\{2, \ldots, d\}$. As a consequence, $\dot{C}_1$ cannot be
extended continuously to the set $\{0\} \times([0, 1]^{d-1} \setminus
(0, 1]^{d-1})$.

In the bivariate case, Condition~\ref{cdiffCK} can be verified by
direct calculation, provided the correlation parameter $\rho$ satisfies
$|\rho| < 1$.
\end{example}

\begin{example}[(Archimedean copulas)]
\label{exarchimedean}
Let $C$ be a $d$-variate Archimedean copula; that is,
\[
C(u) = \phi^{-1}  \bigl( \phi(u_1) + \cdots+ \phi(u_d)  \bigr), \qquad  u
\in[0, 1]^d,
\]
where the generator $\phi\dvtx [0, 1] \to[0, \infty]$ is convex,
decreasing, finite on $(0, 1]$, and vanishes at~$1$, whereas $\phi^{-1}
\dvtx [0, \infty) \to[0, 1]$ is its generalized inverse, $\phi^{-1}(x) =
\inf\{ u \in[0, 1] \dvt \phi(u) \le x \}$; in fact, if $d \ge3$, more
conditions on $\phi$ are required for $C$ to be a copula; see McNeil and Ne{\v{s}}lehov{\'a}~\cite
{mcneilneslehova2009}.

Suppose $\phi$ is continuously differentiable on $(0, 1]$ and $\phi
'(0+) = -\infty$. Then the first-order partial derivatives of $C$ are
given by
\[
\dot{C}_j(u) = \frac{\phi'(u_j)}{\phi'(C(u))},\qquad   u \in[0, 1]^d,
  0 < u_j < 1.
\]
If $u_i = 0$ for some $i \ne j$, then $C(u) = 0$ and $\phi'(C(u)) = -
\infty$, so indeed $\dot{C}_j(u) = 0$. We find that Condition~\ref
{cdiffC} is verified, so Propositions~\ref{pempproc} and~\ref{pmclt} apply.

In contrast, $\dot{C}_j$ may easily fail to be continuous at some
boundary points. For instance, if $\phi'(1) = 0$, then $\dot{C}_j$
cannot be extended continuously at $(1, \ldots, 1)$. Or if $\phi^{-1}$
is long-tailed, that is, if $\lim_{x \to\infty} \phi^{-1}(x+y) / \phi
^{-1}(x) = 1$ for all $y \in\reals$, then $\lim_{u_1 \downarrow0}
C(u_1, u_{-1})/u_1 = 1$ for all $u_{-1} \in(0, 1]^{d-1}$, whereas $\dot
{C}_1(u) = 0$ as soon as $u_j = 0$ for some $j \in\{2, \ldots, d\}$;
it follows that $\dot{C}_1$ cannot be extended continuously to the set
$\{0\} \times([0, 1]^{d-1} \setminus(0, 1]^{d-1})$.
\end{example}

\begin{example}[(Extreme-value copulas)]
\label{exevc}
Let $C$ be a $d$-variate extreme-value copula; that is,
\[
C(u) = \exp ( - \ell( - \log u_1, \ldots, - \log u_d )  ),
 \qquad u \in(0, 1]^d,
\]
where the tail dependence function $\ell\dvtx [0, \infty)^d \to[0, \infty
)$ verifies
\[
\ell(x) = \int_{\Delta_{d-1}} \max_{j \in\{1, \ldots, d\}} (w_j x_j)
  H(\mathrm{d}w),   \qquad x \in[0, \infty)^d,
\]
with $H$ a non-negative Borel measure (called spectral measure) on the
unit simplex $\Delta_{d-1} = \{ w \in[0, 1]^d \dvt w_1 + \cdots+ w_d = 1
\}$ satisfying the $d$ constraints $\int w_j   H(\mathrm{d}w) = 1$ for all $j
\in\{1, \ldots, d\}$; see, for instance, Leadbetter and Rootz{\'e}n \cite
{leadbetterrootzen1988} or Pickands \cite{pickands1989}. It can be verified
that $\ell$ is convex, is homogeneous of order $1$, and that $\max(x_1,
\ldots, x_d) \le\ell(x) \le x_1 + \cdots+ x_d$ for all $x \in[0,
\infty)^d$.

Suppose that the following holds:
%
\begin{equation}
\begin{tabular}{@{}p{330pt}}
\label{cdiffell}
For every $j \in\{1, \ldots, d\}$, the first-order partial derivative
$\dot{\ell}_j$ of $\ell$ with respect to $x_j$ exists and is continuous
on the set $\{ x \in[0, \infty)^d \dvt x_j > 0 \}$.
\end{tabular}
\end{equation}
Then the first-order partial derivative of $C$ in $u$ with respect to
$u_j$ exists and is continuous on the set $\{ u \in[0, 1]^d \dvt 0 < u_j
< 1 \}$. Indeed, for $u \in[0, 1]^d$ such that $0 < u_j < 1$, we have
\[
\dot{C}_j(u) =
\cases{
 \displaystyle\frac{C(u)}{u_j}   \dot{\ell}_j(-\log u_1, \ldots,
-\log u_d), & \quad$\mbox{if $u_i > 0$ for all $i$,} $\vspace*{2pt}\cr
0, & \quad$\mbox{if $u_i = 0$ for some $i \ne j$.}$}
\]
The properties of $\ell$ imply that $0 \le\dot{\ell}_j \le1$ for all
$j \in\{1, \ldots, d\}$. Therefore, if $u_i \downarrow0$ for some $i
\ne j$, then $\dot{C}_j(u) \to0$, as required. Hence if \eqref
{cdiffell} is verified, Condition~\ref{cdiffC} is verified as well
and Propositions~\ref{pempproc} and~\ref{pmclt} apply.

Let us consider the bivariate case in somewhat more detail. The
function $A \dvtx [0, 1] \to[1/2, 1] \dvtx t \mapsto A(t) = \ell(1-t, t)$ is
called the Pickands dependence function of $C$. It is convex and
satisfies $\max(t, 1-t) \le A(t) \le1$ for all $t \in[0, 1]$. By
homogeneity of the function~$\ell$, we have $\ell(x, y) = (x+y)
A(\frac{y}{x+y})$ for $(x, y) \in[0, \infty)^2 \setminus\{(0, 0)\}$.
If $A$ is continuously differentiable on $(0, 1)$, then \eqref
{cdiffell} holds, and Condition~\ref{cdiffC} is verified.
Nevertheless, if $A(1/2) < 1$, which is always true except in case of
independence ($A \equiv1$), the upper tail dependence coefficient $2\{
1 - A(1/2)\}$ is positive so that the first-order partial derivatives
fail to be continuous at the point $(1, 1)$; see Example~\ref
{extaildep}. One can also see that $\dot{C}_1$ will not admit a
continuous extension in the neighborhood of the point $(0, 0)$ in case
$A'(0) = -1$.

We will now verify Condition~\ref{cdiffCK} under the following
additional assumption:
%
\begin{equation}
\begin{tabular}{@{}p{244pt}}
\label{cdiffA}
The function $A$ is twice continuously differentiable on $(0, 1)$ and
$M = \sup_{0 < t < 1} \{ t   (1-t)   A''(t) \} < \infty$.
\end{tabular}
\end{equation}
In combination with Proposition~\ref{pempprocrate}, this will justify
the use of the Stute--Tsukahara almost sure rate~\eqref{estute} in the
proof of Theorem~3.2 in Genest and Segers \cite{genestsegers2009}; in particular, see
their equation~(B.3). Note that the weight function $t   (1-t)$ in the
supremum in \eqref{cdiffA} is not unimportant: for the Gumbel
extreme-value copula having dependence function $A(t) = \{ t^{1/\theta}
+ (1-t)^{1/\theta} \}^\theta$ with parameter $\theta\in(0, 1]$, it
holds that $A''(t) \to\infty$ as $t \to0$ or $t \to1$ provided $1/2
< \theta< 1$, whereas condition~\eqref{cdiffA} is verified for all
$\theta\in(0, 1]$.

The copula density at the point $(u, v) \in(0, 1)^2$ is given by
\[
\ddot{C}_{12}(u, v) = \frac{C(u, v)}{uv} \biggl( \mu(t)   \nu(t) -
\frac{t   (1-t)   A''(t)}{\log(uv)} \biggr),
\]
where
\[
t = \frac{\log(v)}{\log(uv)} \in(0, 1),\qquad \mu(t) = A(t) - t
A'(t),  \nu(t) = A(t) + (1-t)   A'(t).
\]
Note that if $A''(1/2) > 0$, then $\ddot{C}_{12}(w, w) \to\infty$ as
$w \uparrow1$. The properties of $A$ imply $0 \le\mu(t) \le1$ and $0
\le\nu(t) \le1$. From $-\log(x) \ge1-x$, it follows that $-1/\log
(uv) \le\min\{ 1/(1-u), 1/(1-v) \}$ for $(u, v) \in(0, 1)^2$. Since
$C(u, v) \le\min(u, v)$ and since $\min(a, b)   \min(c, d) \le\min
\{ (ac), (bd) \}$ for positive numbers $a, b, c, d$, we find
\begin{eqnarray*}
0 &\le&\ddot{C}_{12}(u, v)
\le\frac{\min(u, v)}{uv} \biggl\{ 1 + M   \min\biggl( \frac
{1}{1-u}, \frac{1}{1-v} \biggr) \biggr\} \\
&\le&(1 + M)   \min\biggl( \frac{1}{u(1-u)} , \frac{1}{v(1-v)} \biggr).
\end{eqnarray*}
Similarly, for $(u, v) \in(0, 1) \times[0, 1]$,
\[
\ddot{C}_{11}(u, v) =
\cases{
 \displaystyle\frac{C(u,v)}{u^2}   \biggl( - \mu(t)    \bigl( 1 - \mu
(t)  \bigr) + \displaystyle\frac{t^2   (1-t)   A''(t)}{\log(u)} \biggr), &\quad$\mbox{if
$0 < v < 1$},$\vspace*{2pt}\cr
0, & \quad$\mbox{if $v \in\{0, 1\}$}.$}
\]
Continuity at the boundary $v = 0$ follows from the fact that $C(u, v)
\to0$ as $v \to0$; continuity at the boundary $v = 1$ follows from
the fact that $t \to0$ and $\mu(t) \to0$ as $v \to1$. Since $- \log
(u) \le(1-u)/u$, we find, as required,
\[
0 \le- \ddot{C}_{11}(u, v)
\le\frac{(1 + M)}{u   (1-u)},\qquad
 (u, v) \in(0, 1) \times[0, 1].
\]
\end{example}

\begin{example}[(If everything fails\ldots)]
Sometimes, even Condition~\ref{cdiffC} does not hold: think, for
instance, of the Fr\'echet lower and upper bounds, $C(u, v) = \max(u +
v - 1, 0)$ and $C(u, v) = \min(u, v)$,\vadjust{\goodbreak} and of the checkerboard copula
with Lebesgue density $c = 2   \ind_{[0, 1/2]^2 \cup[1/2, 1]^2}$. In
these cases, the candidate limiting process $\CC$ has discontinuous
trajectories, and the empirical copula process does not converge weakly
in the topology of uniform convergence.

One may then wonder if weak convergence of the empirical copula process
still holds in, for instance, a Skorohod-type topology on the space of
c\`{a}dl\`{a}g functions on $[0, 1]^2$. Such a result would be useful to
derive, for instance, the asymptotic distribution of certain
functionals of the empirical copula process, for example, suprema or
integrals such as appearing in certain test statistics.
\end{example}

\begin{appendix}\label{app}

\section*{Appendix: Multivariate oscillation modulus}

Let $C$ be any $d$-variate copula and let $U_1, U_2, \ldots$ be an
i.i.d. sequence of random vectors with common cumulative distribution
function $C$. Let $\alpha_n$ be the multivariate empirical process in
equation~\eqref{ealphan}. Consider the oscillation modulus defined by
\setcounter{equation}{0}
\begin{equation}
\label{eoscillation}
M_n(a) = \sup\{ | \alpha_n(u) - \alpha_n(v) | \dvt \mbox{$u, v \in[0,
1]^d$, $|u_j-v_j| \le a$ for all $j$} \}
\end{equation}
for $a \in[0, \infty)$. Define the function $\psi\dvtx [-1,\infty) \to
(0, \infty)$ by
%
\begin{equation}
\label{epsi}
\psi(x) = 2x^{-2} \{ (1+x) \log(1+x) - x \}, \qquad  x \in(-1, 0) \cup
(0, \infty),
\end{equation}
together with $\psi(-1) = 2$ and $\psi(0) = 1$. Note that $\psi$ is
decreasing and continuous.

\setcounter{proposition}{0}
\begin{proposition}[(John H. J. Einmahl, Hideatsu Tsukahara)]
\label{poscillation}
Let $C$ be any $d$-variate copula. There exist constants $K_1$ and
$K_2$, depending only on $d$, such that
\[
\Pr\{ M_n(a) \ge\lambda\} \le\frac{K_1}{a}   \exp\biggl\{ - \frac
{K_2 \lambda^2}{a}   \psi\biggl( \frac{\lambda}{\sqrt{n} a} \biggr)
\biggr\}
\]
for all $a \in(0, 1/2]$ and all $\lambda\in[0, \infty)$.
\end{proposition}

\begin{pf}
In Einmahl \cite{einmahl1987}, Inequality~5.3, page 73, the same bound is
proved in case $C$ is the independence copula and for $a > 0$ such that
$1/a$ is integer. As noted by Tsukahara, in a private communication,
the only property of the joint distribution that is used in the proof
is that the margins be uniform on the interval $(0, 1)$: Inequality~2.5
in Einmahl \cite{einmahl1987}, page~12, holds for any distribution on the unit
hypercube and equation~(5.19) on page~72 only involves the margins. As
a consequence, Inequality~5.3 in Einmahl \cite{einmahl1987} continues to hold
for any copula $C$. Moreover, the assumption that $1/a$ be integer is
easy to get rid of.
\end{pf}
\end{appendix}

\section*{Acknowledgments}

Input from an exceptionally large number of sources has shaped the
paper in its present form. The author is indebted to the following persons:
\textit{Axel B\"ucher}, for careful reading and for mentioning the
possibility of an alternative proof of Proposition~\ref{pempproc} via
the functional delta method as in B\"{u}cher  \cite{buecher2011};
\textit{John H.J. Einmahl}, for pointing out the almost sure bound in
equation~\eqref{ebound1} on the tail empirical process and the
resulting reinforcement of the conclusion of Proposition~\ref
{pempprocrate} with respect to an earlier version of the paper;
\textit{Gordon Gudendorf}, for fruitful discussions on Condition~\ref
{cdiffCK} in the context of extreme-value copulas;
\textit{Ivan Kojadinovic}, for meticulous proofreading resulting in a~substantially lower error rate and for numerous suggestions leading to
improvements in the exposition;
\textit{Hideatsu Tsukahara}, for sharing the correction of the final part
of the proof of Theorem~4.1 in Tsukahara \cite{tsukahara2000}, summarized here
in the proof of Proposition~\ref{pempprocrate} and in Proposition~\ref
{poscillation};
\textit{the referees, the Associate Editor, and the Editor}, for timely
and constructive comments.
Funding was provided by IAP research network Grant P6/03 of the Belgian
government (Belgian Science Policy)
and by ``Projet d'actions de recherche concert\'ees'' number 07/12/002
of the Communaut\'e fran\c{c}aise de
Belgique, granted by the Acad\'emie universitaire de Louvain.


%

\printhistory

\end{document}